\documentstyle[twoside,12pt]{article}
\newtheorem{prop}{Proposition}[section]

\newtheorem{rem}[prop]{Remark}
\newtheorem{coro}[prop]{Corollary}

\title{On Betti numbers and Chern classes 
of varieties with trivial odd cohomology groups}

\author{Lev A. Borisov \\
Mathematical Sciences\\
Research Institute \\
Berkeley, CA 94720 \\
e-mail: lborisov@msri.org}

\begin{document}

\date{}

\maketitle

\begin{abstract}
It was observed in a very recent preprint (hep-th/9703086) 
by T. Eguchi, K. Hori, and Ch.-Sh. Xiong that a curious identity 
between Betti numbers and Chern classes holds
for many examples of Fano varieties. The goal of this paper is
to prove that for varieties with trivial odd cohomology groups
this identity is equivalent to having zero Hodge numbers $h^{p,q}$
for $p\neq q$.
\end{abstract}

\section{Introduction}

\noindent

Let $X$ be a smooth complex projective variety of dimension $n$ whose odd cohomology
groups $H^{2k+1}({\bf C})$ are zero. It was noticed in \cite{tekhcx} that 
a curious identity
$$\frac 14 \sum_k h^{2k} (k-\frac{n-1}2)(1-k+\frac{n-1}2)=
\frac1{24}\left(\frac{3-n}2\chi(X)-\int_Xc_1(X)\wedge c_{n-1}(X)\right)$$
holds in many examples. The authors of \cite{tekhcx} were primarily
concerned with Fano varieties $X$ and the above identity was a prerequisite
for the conjectural construction of Virasoro operators that control
quantum cohomology of $X$.
The above condition could be rewritten as
$$\sum_{k=0}^n h^{2k} (k-\frac n2)^2=\frac16 c_1c_{n-1}+\frac n{12} c_n$$
where $c_l$ is the $l$-th Chern class of $X$.
The goal of this paper is to prove the following result.

\begin{prop}
If $X$ is a smooth complex projective variety of dimension $n$ with $H^{odd}(X,{\bf C})=0$
then
$$\sum_{k=0}^n h^{2k} (k-\frac n2)^2\leq\frac16 c_1c_{n-1}+\frac n{12} c_n$$
and equality holds if and only if
$$h^{p,q}=0 ~{\rm for~all}~p\neq q.$$
\end{prop}
\label{main}

In addition, we give an application of this result to the combinatorics
of reflexive polytopes that describe smooth toric Fano varieties.

After this preprint was submitted to the archive, the author was informed by 
Anatoly Libgober that the proof of the crucial Proposition 2.2 is contained
in \cite{liwo}. In addition, Victor Batyrev presented to the author
a combinatorial proof of Corollary 2.3 for arbitrary smooth toric 
varieties.


\section{Proof of the main result}

\noindent

Let $X$ be a smooth complex projective variety of dimension $n$.
We introduce the E-polynomial
$$E(u,v)=\sum_{p,q}(-1)^{p+q}h^{p,q}u^pv^q$$
where $h^{p,q}={\rm dim}H^q(X,\Lambda^p T^*X)$ are Hodge numbers of $X$.
We also introduce
$$\chi_p=(-1)^p\chi(\Lambda^p T^*X)=\sum_q(-1)^{p+q}h^{p,q}$$
and the polynomial 
$$\hat E(t)=\sum_p\chi_p t^p=E(t,1).$$

\begin{rem}
If $h^{p,q}=0$ for all $p\neq q$ then $\chi_p=h^{p,p}$.
\end{rem}
\label{hpp}

\begin{prop}
In the above notations
$$\sum_{p=0}^n\chi_p(p-\frac n2)^2=\frac16c_1c_{n-1}+\frac n{12}c_n.$$
\end{prop}
\label{chiequal}

{\em Proof.} This result comes as an easy application of Hirzebruch-Riemann-Roch
theorem \cite{hirz}. We start by rewriting the left hand side in terms of the 
polynomial $\hat E$.
$$\sum_p\chi_p(p-\frac n2)^2=\sum_p\chi_p p(p-1) +
\sum_p \chi_p(1-n)(p-\frac n2)+(\sum_p \chi_p)(\frac n2-\frac {n^2}4)$$
$$=\frac {d^2}{dt^2}\hat E(t)|_{t=1}+\hat E(1)(\frac n2-\frac {n^2}4).$$
We have used $\chi_p=\chi_{n-p}$ to get rid of the second sum.
By Hirzebruch-Riemann-Roch theorem,
$$\hat E(t)=\sum_p t^p(-1)^p\chi(\Lambda^pT^*X)=
\int_X Td(X)\sum_p(-t)^p ch(\Lambda^pT^*X).$$
As usual, we introduce Chern roots $\alpha_i$ such that
$c(TX)(w)=\prod_i(1+\alpha_iw)$, see for example \cite{fult}.
Then 
$$\hat E(t)=\int_X (\prod_i\frac{\alpha_i}{1-{\rm e}^{-\alpha_i}})
\sum_p(-t)^p\sum_{i_1<...<i_p}{\rm e}^{-\alpha_{i_1}-...-\alpha_{i_p}}$$
$$=\int_X \prod_i \alpha_i(1+(1-t)\frac{e^{-\alpha_i}}{1-{\rm e}^{-\alpha_i}}).$$
This shows, of course, that $\hat E(1)=\chi(X)=c_n$.
Besides we can calculate
$$
\frac {d^2}{dt^2}\hat E(t)|_{t=1}=
2\sum_{i<j}\int_X (\prod_{k\neq i,j} \alpha_k)
(1-\frac12\alpha_i+\frac 1{12}\alpha_i^2)(1-\frac12\alpha_j+\frac 1{12}\alpha_j^2)
$$
$$=\frac16 c_1c_{n-1}+(\frac {n^2}4-\frac {5n}{12})c_n$$
and the rest is straightforward.
\hfill $\Box$

We combine Proposition 2.2 with Remark 2.1 to get the following
corollary.

\begin{coro}
If $h^{p,q}=0$ for all $p\neq q$ then
$$\sum_{p=0}^n h^{2p}(p-\frac n2)^2=\frac16c_1c_{n-1}+\frac n{12}c_n.$$
\end{coro}

This corollary gives a sufficient condition for the identity of \cite{tekhcx}.
Now we will see that this condition is also necessary, that is we will prove
Proposition 1.1.

{\em Proof of Proposition 1.1.}
Because of $h^{odd}=0$, we have $\chi_p=\sum_q h^{p,q}$ and  
$$\sum_{p=0}^n h^{2p}(p-\frac n2)^2=
\sum_{p,q}h^{p,q}(\frac{p+q}2-\frac n2)^2 $$
$$=\sum_p (\sum_q h^{p,q})(p-\frac n2)^2+
\sum_{p,q} h^{p,q}(\frac {q-p}2)(\frac {3p+q}2-n)$$
$$
=\sum_p \chi_p (p-\frac n2)^2-
\sum_{p,q}h^{p,q}(\frac {q-p}2)^2+
\sum_{p,q}h^{p,q}(\frac {q-p}2)(p+q-n).$$
Because of $h^{p,q}=h^{q,p}$, the last sum is zero. Together
with the result of Proposition 2.2, we get
$$\sum_{p=0}^n h^{2p}(p-\frac n2)^2=
\frac 16 c_1c_{n-1}+\frac 1{12} c_n-\sum_{p,q}h^{p,q}(\frac {q-p}2)^2$$
which proves the proposition. \hfill ${\Box}$

\section{Application to toric Fano varieties}
The goal of this section is to give a combinatorial equivalent of
Corollary 2.3  for the case of smooth toric Fano varieties.

Recall (see for example \cite{baty}) that a smooth toric Fano variety
can be defined in terms of the polytope $\Delta \in M$ that supports 
the sections of the anticanonical line bundle. 

To calculate the left hand side, notice that $X$ is
a disjoint union of algebraic tori $T_\theta$. Using the additivity
of the cohomology with compact support, we get
$$\hat E(t)=\sum_{\theta\subseteq \Delta}(t-1)^{{\rm dim}(\theta)}$$
and 
$$\frac {d^2}{dt^2}\hat E(t)|_{t=1}=2\#\{\theta\subseteq \Delta, {\rm dim}(\theta)=2\}.$$
Analogously, 
$$\hat E(1)=\#\{\theta\subseteq \Delta, {\rm dim}(\theta)=0\}.$$
To calculate the right hand side of the identity, notice that
$$c(TX)(w)=\prod_{\theta, {\rm dim}(\theta)=n-1} (1+D_\theta w),$$
where $D_\theta$ is the closure of the strata that corresponds to $\theta$,
and that $\sum D_\theta$ is a divisor with normal 
crossings. This gives 
$$c_1=\sum_{\theta, {\rm dim}(\theta)=n-1} D_\theta$$
$$c_{n-1}=\sum_{\theta, {\rm dim}(\theta)=1} l_\theta$$
where in the second identity $l_\theta$ is the closure of $T_\theta$.
One can show that $c_1 l_\theta$ equals the number of interior
points of $\theta$ plus one. Now an easy calculation shows that
the identity of Corollary 2.3 becomes
$$
\#\{\theta\subseteq \Delta, {\rm dim}(\theta)=2\}=
\frac1{12}\sum_{\theta,{\rm dim}(\theta)=1}
\#\{P\in M,~P\in {\rm interior}(\theta)\}$$
$$+
(\frac {n^2}8-\frac n6)\#\{\theta\subseteq \Delta, {\rm dim}(\theta)=0\}.
$$

\end{document}